\documentclass[12pt]{article}
\usepackage{amsmath,amsfonts,amssymb,amsthm}
\begin{document}

\newcommand{\1}{{{\bf 1}}}
\newcommand{\id}{{\rm id}}
\newcommand{\Hom}{{\rm Hom}\,}
\newcommand{\End}{{\rm End}\,}
\newcommand{\Res}{{\rm Res}\,}
\newcommand{\Image}{{\rm Im}\,}
\newcommand{\Ind}{{\rm Ind}\,}
\newcommand{\Aut}{{\rm Aut}\,}
\newcommand{\Ker}{{\rm Ker}\,}
\newcommand{\gr}{{\rm gr}}
\newcommand{\Der}{{\rm Der}\,}
\newcommand{\Span}{{\rm Span}}
\newcommand{\tr}{{\rm tr}}

\newcommand{\Z}{\mathbb{Z}}
\newcommand{\Q}{\mathbb{Q}}
\newcommand{\C}{\mathbb{C}}
\newcommand{\N}{\mathbb{N}}

\newcommand{\g}{\mathfrak{g}}
\newcommand{\h}{\mathfrak{h}}
\newcommand{\wt}{{\rm wt}\;}
\newcommand{\CR}{\mathcal{R}}
\newcommand{\D}{\mathcal{D}}
\newcommand{\E}{\mathcal{E}}
\newcommand{\Lie}{\mathcal{L}}
\newcommand{\z}{\bf{z}}
\newcommand{\bflam}{\bf{\lambda}}

\newcommand{\levelk}{\bf{k}}

\def \<{\langle} 
\def \>{\rangle}
\def \be{\begin{equation}\label}
\def \ee{\end{equation}}
\def \bex{\begin{exa}\label}
\def \eex{\end{exa}}
\def \bl{\begin{lem}\label}
\def \el{\end{lem}}
\def \bt{\begin{thm}\label}
\def \et{\end{thm}}
\def \bp{\begin{prop}\label}
\def \ep{\end{prop}}
\def \br{\begin{rem}\label}
\def \er{\end{rem}}
\def \bc{\begin{coro}\label}
\def \ec{\end{coro}}
\def \bd{\begin{de}\label}
\def \ed{\end{de}}

\def \bproof{\noindent \textbf{Proof:\ \ }}
\def \eproof{$\blacksquare$}

\newtheorem{thm}{Theorem}[section]
\newtheorem{prop}[thm]{Proposition}
\newtheorem{coro}[thm]{Corollary}
\newtheorem{conj}[thm]{Conjecture}
\newtheorem{exa}[thm]{Example}
\newtheorem{lem}[thm]{Lemma}
\newtheorem{rem}[thm]{Remark}
\newtheorem{de}[thm]{Definition}
\newtheorem{hy}[thm]{Hypothesis}
\makeatletter
\@addtoreset{equation}{section}
\def\theequation{\thesection.\arabic{equation}}
\makeatother
\makeatletter


\title{A recurrence relation
for characters of highest weight integrable modules for
affine Lie algebras}
\author{William J. Cook\footnotemark[1], Haisheng Li\footnotemark[2], 
and Kailash C. Misra\footnotemark[3]
\\
{\small \footnotemark[1] Department of Mathematics, 
Rutgers University, Piscataway, NJ 08854}\\
{\small \footnotemark[2] Department of Mathematical Sciences, 
Rutgers University, Camden, NJ 08102}\\
{\small \footnotemark[3] Department of Mathematics, North Carolina State University, Raleigh, NC 27695}} 

\renewcommand{\thefootnote}{\fnsymbol{footnote}}
\footnotetext[0]{Cook acknowledges partial support by NSA Grant MDA
  904-02-1-0072 and NSF Grant CCR-0096842. Misra acknowledges the
  partial support by the NSA grants MDA904-02-1-0072 and H98230-06-1-0025. 
Li acknowledges partial support by the NSA grant H98230-05-1-0018.}

\maketitle 
 
\begin{abstract}
Using certain results for the vertex operator algebras associated with
affine Lie algebras we obtain recurrence relations
for the characters of integrable highest weight irreducible
modules for an affine Lie algebra.  As an application 
we show that in the simply-laced level $1$ case, 
these recurrence relations give the known characters, 
whose principal specializations naturally give rise to 
some multisum Macdonald identities.
\end{abstract}

\section{Introduction} 
Interactions between representations of affine Lie algebras
and combinatorial identities have been well known for
many years. In 1972, Macdonald \cite{MacD} obtained a family of
identities related to the affine root systems (see also \cite{dyson}). This led to
further investigations (see \cite{k1974}, \cite{M}, \cite{L1}, \cite{kvstrange}, \cite{L2})
on combinatorial identities and affine Lie algebras. 
In 1978, Lepowsky and Milne \cite{LM} showed that the
product sides of the famous Rogers-Ramanujan identities appear in the
principal specialization of the characters of certain irreducible
representations of the simplest affine Lie algebra
$\widehat{sl_{2}}$. Further investigation of this observation led Lepowsky
and Wilson \cite{LW1} to give a vertex operator construction of the
basic representations of $\widehat{sl_{2}}$, which was used later by them
to give a vertex-operator theoretic proof of the Rogers-Ramanujan identities
\cite{LW3}. In \cite{LW2}, Lepowsky and Wilson introduced 
a new family of algebras, which they called
$Z$-algebras, which provided a new tool for the investigation of the
connections between highest weight integrable representations 
for affine Lie algebras  and
combinatorial identities. These investigations subsequently led to new
combinatorial identities and new algebraic structures such as vertex
(operator) algebras (\cite{B}, \cite{flm}).

It is well known (\cite{B}, \cite{flm}, \cite{fz}, \cite{dl}) that the irreducible
integrable highest weight module $L(k,0)$ $(=L(k\Lambda_0))$ of level
$k>0$ for an affine Lie algebra $\hat{\g}$ has a natural structure of a
vertex operator algebra.  For many years, vertex operator methods and vertex operator algebras
have been successfully used as new tools in the study of representations
for affine Lie algebras (cf. \cite{c}, \cite{LP}, \cite{LW4}, \cite{LW5}, \cite{mpaim} -
\cite{mpbook}, \cite{Mi1}, \cite{Mi2}, \cite{primc9805}, \cite{primccombiden}).  
A very exciting example is the recent work (\cite{clm1}, \cite{clm2}), in which by using the
vertex operator algebras $L(k,0)$ associated with $\widehat{sl_{2}}$ 
and certain intertwining operators, Capparelli, Lepowsky and Milas
recovered the celebrated Rogers-Ramanujan recursions and Rogers-Selberg
recursions in relation to certain subspaces of the integrable highest weight modules.

This paper is mostly motivated by 
the work \cite{clm1} and \cite{clm2}. In this paper, using certain results of \cite{li-affine} for the vertex operator algebra $L(k,0)$  we obtain certain recurrence relations for the full characters of the integrable highest weight
$\hat{\g}$-modules. 
Although motivated by \cite{clm1} and \cite{clm2}, our methods for obtaining recursions and the recursions themselves
are quite different. Those that we derive are related to the invariance of characters under the
translations in the affine Weyl group. We are essentially interpreting this invariance using vertex-operator 
algebraic methods.
As an application for the case that $\g$ is simply-laced and $k=1$, we use these recurrence relations to recover known multisum expressions 
for the characters of these affine Lie algebra modules. In particular,  taking the principal specializations of these
characters and equating them with the known product formulas \cite{L3}, we easily obtain some combinatorial identities which had been proved by Macdonald in \cite{MacD}. In the case of the affine Lie algebra $A_{\ell}^{(1)}$ it has been
pointed out to us by Ole Warnaar that the multisum expression for the principally specialized character is the generating function for $q$-cores \cite{GKS}. For the affine Lie algebra $\widehat{sl_{2}}$ this amounts to a well known identity due to Gauss (cf. \cite{A}, \cite{fl}, \cite{kvstrange}). 

The authors would like to thank James Lepowsky and Antun Milas for their insightful comments and suggestions which
have greatly improved the presentation of this work.  

\section{Recurrence relations for characters of irreducible $L(k,0)$-modules}

In this section, we first review some basic notations and results about
(affine) Lie algebras and vertex operator algebras, 
and then using some results in
\cite{li-affine} on vertex operator algebra representations we derive
recurrence relations for the characters of integrable highest weight
irreducible modules for the affine Lie algebra $\hat{\g}$.

First we review some basic facts
about finite-dimensional simple Lie algebras (cf. \cite{hum})
and (untwisted) affine Lie algebras (cf. \cite{kacbook1}).
Let $\g$ be a finite-dimensional simple Lie algebra of rank $\ell$,
let $\h$ be a fixed Cartan subalgebra and let $\<\cdot,\cdot\>$ 
be the normalized Killing form such that $\<\alpha,\alpha\>=2$
for long roots $\alpha$. Also, let $\Delta$ denote the set of roots of
$\g$ and fix a set $\{\alpha_1,\dots,\alpha_{\ell}\}$ of simple roots. 

Let $\lambda_i$ for $i=1,\dots,\ell$ be the fundamental weights of $\g$. 
Denote by $Q$ and $P$ the 
root lattice and the weight lattice of $\g$:
\begin{eqnarray}
Q&=& \Z \alpha_1 + \dots + \Z \alpha_{\ell}, \\
P&=& \Z \lambda_1 + \dots + \Z \lambda_{\ell}.
\end{eqnarray}

For $\alpha=m_{1}\alpha_{1}+\dots +m_{\ell}\alpha_{\ell}\in Q$, set
\begin{eqnarray}
{\rm ht}(\alpha)=m_{1}+\cdots +m_{\ell}.
\end{eqnarray}

Let $\{E_i,F_i,H_i \,|\, 1 \leq i \leq \ell \}$ 
be a set of Chevalley generators of $\g$.
We have $\alpha_j(H_i) = a_{ij}$ where
$C=(a_{ij})_{i,j=1}^\ell$ is the Cartan matrix of $\g$.  
Let $H^{(i)}$ for $i=1,\dots, \ell$ be the fundamental
co-weights, i.e., $H^{(j)}\in \h$ such that
$\alpha_i(H^{(j)})=\delta_{i,j}$ for $i=1,\dots,\ell$. The fundamental
co-weights form another basis of $\h$. For $1\le i,j\le \ell$, we have
\begin{eqnarray}
\< H_i, H^{(j)}\> = \delta_{i,j} \frac{2}{\<\alpha_i, \alpha_i\>}.
\end{eqnarray}

\noindent
Denote by $Q^{\vee}$ and $P^{\vee}$ the 
co-root lattice and the co-weight lattice of $\g$. That is,
\begin{eqnarray}
Q^{\vee} &=& \Z H_1 + \cdots +\Z H_{\ell},\\
P^{\vee} &=& \Z H^{(1)} + \cdots +\Z H^{(\ell)}.
\end{eqnarray}

\noindent
Set
\begin{eqnarray}
\rho^{\vee} = \sum_{i=1}^{\ell} H^{(i)} \in \h.
\end{eqnarray}

\noindent
We have
\begin{equation}
\alpha_j(\rho^{\vee}) = \sum_{i=1}^{\ell}\alpha_j (H^{(i)}) 
                      = \sum_{i=1}^{\ell} \delta_{i,j} = 1
\end{equation}

\noindent
for all $1 \leq j \leq \ell$. Consequently,
\begin{eqnarray}
\alpha(\rho^{\vee})={\rm ht}(\alpha)
\;\;\;\mbox{ for }\alpha\in \Delta.
\end{eqnarray}

Recall the (untwisted) affine Lie algebra
\begin{eqnarray}
\hat{\g} = \g \otimes \C [t,t^{-1}]\oplus \C c,
\end{eqnarray}

\noindent
where for $a,b \in \g$, $m,n \in \Z$,
\begin{eqnarray}
& &[a\otimes t^m, b\otimes t^n] 
 = [a,b] \otimes t^{m+n} + m \<a,b\>\delta_{m+n,0} c,\\
& &\hspace{2cm} \ \ \ \ \  [\hat{\g},c]=0.
\end{eqnarray}

For $a\in \g$, $n \in \Z$, denote by $a(n)$ the 
corresponding operator of $a\otimes t^{n}$ on $\hat{\g}$-modules.
Let $\theta$ be the highest long root of $\g$. 
Choose (non-zero) vectors $E_{\theta} \in \g_{\theta}$ and $F_{\theta}
\in \g_{-\theta}$ such that $\< E_{\theta}, F_{\theta}\>=1$.  Set
$H_{\theta} = [ E_{\theta},F_{\theta}]$.
Set $e_0 = F_{\theta} \otimes t$, $f_0 = E_{\theta} \otimes t^{-1}$,
and $h_0 = [e_0,f_0]$.  Also, let $e_i = E_i \otimes 1$, $f_i = F_i
\otimes 1$, and $h_i = H_i \otimes 1$ for $i=1,\dots,\ell$. {}From
\cite{kacbook1}, $\{ e_i, f_i, h_i \,|\, 0 \leq i \leq \ell\}$
is a set of Chevalley generators of $\hat{\g}$.

Let $\tilde{\g} = \hat{\g} \oplus \C d$ be the extended affine Lie
algebra (semi-direct product Lie algebra), where $d$ is the derivation
of $\hat{\g}$ acting on $\g \otimes \C [t,t^{-1}]$ as $1 \otimes t
\frac{d}{dt}$ and acting on $\C c$ as $0$. The subalgebra $\h\oplus\C
c\oplus\C d$ is a Cartan subalgebra of $\tilde{\g}$. The
\emph{homogeneous} grading and the \emph{principal} grading on $\hat\g$ (and
on $\tilde{\g}$) are given by the derivations $d_H$ and $d_P$, respectively, where
\begin{eqnarray}
& & d_H(e_0) = e_0,\;\;\;\; d_H(f_0)=-f_0,\;\;\;\; d_H(c)=0, \\
& & d_H(e_i) = 0 = d_H(f_i) 
   \; \hbox{ for } \; 1 \leq i \leq \ell \notag
\end{eqnarray}

\noindent
and
\begin{eqnarray}
& & d_P(e_i)=e_i,\;\;\;\; d_P(f_i)=-f_i 
     \; \mbox{  for } \; 0 \leq i \leq \ell,\\
& & d_P(c)=0. \notag
\end{eqnarray}

In the homogeneous (resp. principal) grading we have
\begin{equation}
\mathrm{deg}(c)=0, \; \mathrm{deg}(a \otimes t^n) = n \;
\hbox{ (resp. $n \mathrm{ht}(\theta) + \mathrm{ht}(\alpha)$)}  
\end{equation}

\noindent
for all $a \in \g$ (resp. $a \in \g_{\alpha}$, $\alpha \in \Delta$) 
and $n \in \Z$. Observe that $d_H = -d$. We also have the following known relation
between $d_P$ and $d$, which is easy to prove.

\bl{principal-grading} As derivations on $\hat{\g}$, 
\begin{eqnarray}
d_P = (\mathrm{ht}(\theta)+1)d+\mathrm{ad}\;\rho^{\vee}=
-(\mathrm{ht}(\theta)+1)d_H+\mathrm{ad}\; \rho^{\vee}.
\end{eqnarray}

\el

For any $\lambda\in \h^{*},\;k \in \C$, we extend $\lambda$ to a
linear functional $(k,\lambda)$ on $\h \oplus \C c \oplus \C d$ by
defining $(k,\lambda)(c)=k,\; (k,\lambda)(d)=0$.  For $\lambda\in
\h^{*},k \in \C$, we denote by $L(k,\lambda)$ the irreducible highest
weight $\tilde{\g}$-module (of level $k$) with highest weight
$(k,\lambda)$.  
On the irreducible highest weight $\tilde{\g}$-module
$L(k,\lambda)$ for $k \in \C,\; \lambda\in \h^{*}$
we have the homogeneous grading and principal grading,
where we assign the degree of the
highest weight subspace of $L(k,\lambda)$ to be zero.

For $0\le i\le \ell$, set $\Lambda_i
=(1,\lambda_{i})\in (\h \oplus \C c \oplus \C d)^{\ast}$, where
$\lambda_{i}$ for $1\le i\le \ell$ are the fundamental weights of $\g$
and $\lambda_{0}=0$. In terms of these notations we in particular have
$L(k,0)=L(k \Lambda_0)$.

Next we review the basic results on vertex operator algebras
associated with the affine Lie algebra $\hat{\g}$. 
First we recall the definitions from \cite{flm} and \cite{fhl} (cf. \cite{ll-book}).
A {\em vertex operator algebra} is a $\Z$-graded 
vector space $V = \coprod_{n \in \Z} V_{(n)}$ (over $\C$) 
such that  $\dim V_{(n)}<\infty$ for all $n\in \Z$ and
$V_{(n)} = 0$ for $n$ sufficiently small,
which is equipped with a linear map, called the {\em vertex operator map},
\begin{eqnarray}
Y(\cdot,x):V & \rightarrow & (\End V)[[x,x^{-1}]] \notag \\
           v & \mapsto     & Y(v,x) = \sum_{n \in \Z} v_n x^{-n-1},
\end{eqnarray}

\noindent
and equipped with distinguished vectors
$\1 \in V_{(0)}$ (the \emph{vacuum vector}) and $\omega \in V_{(2)}$ (the 
\emph{conformal vector}) such that for $u,v\in V$,
$$u_{n}v=0\;\;\;\mbox{ for $n$ sufficiently large},$$
$$Y(\1,x)v = v,$$                                                                                                
\begin{equation} Y(v,x)\1 \in V[[x]]\;\;\;\mbox{ and }\;
\;\;\lim_{x \rightarrow 0} Y(v,x){\bf 1} = v_{-1} {\bf 1} = v, \label{creationprop} \end{equation}
\begin{eqnarray}
& &x_0^{-1} \delta\left( \frac{x_1 - x_2}{x_0} \right) Y(u,x_1)Y(v,x_2) - 
   x_0^{-1} \delta\left( \frac{x_2 - x_1}{-x_0} \right)
   Y(v,x_2)Y(u,x_1) \nonumber\\
& &\ \ \ \ \ \ 
= x_2^{-1} \delta\left( \frac{x_1 - x_0}{x_2} \right) Y(Y(u,x_0)v,x_2) 
\end{eqnarray}   
(the \emph{Jacobi identity}). 
The following \emph{Virasoro} relation is assumed
\begin{equation} 
[L(m),L(n)] = (m-n)L(m+n) + \frac{m^3-m}{12} \delta_{m+n,0} \, c_V  
\end{equation}

\noindent
for $m,n \in \Z$, where
\begin{equation}
Y(\omega,x) = \sum_{n \in \Z} \omega_n x^{-n-1} 
= \sum_{n \in \Z} L(n) x^{-n-2}
\end{equation}
and where $c_V \in \C$ is called the \emph{central charge} 
(or \emph{rank}) of $V$. 
It is also assumed that
\begin{eqnarray} 
 & &Y(L(-1)v,x) = \frac{d}{dx} Y(v,x)
\;\;\;\mbox{ for }v\in V,\\
& & V_{(n)} = \{ v \in V \,|\, L(0)v = nv \} \;\;\;\mbox{ for }n \in \Z.
\end{eqnarray}

\noindent
The vertex operator algebra $V$ is sometimes denoted by the quadruple 
$(V,Y,\1,\omega)$.

For a vertex operator algebra $V$, a $V$-module is defined (cf. \cite{fhl}, \cite{ll-book}) by using all
the axioms above that make sense, except that a $V$-module $W$ is $\C$-graded
$W=\coprod_{h\in \C}W_{(h)}$, instead of $\Z$-graded, such that
for any $h\in \C$, $W_{(h+n)}=0$ for sufficiently negative integers $n$. 
If $W$ is an indecomposible $V$-module, e.g., an irreducible $V$-module, then
$W=\coprod_{n\in \N}W_{(h+n)}$ for some $h\in \C$. 









For the rest of this paper we assume that $k$ is a positive
integer.  It is well known 
(see \cite{fz}, \cite{dl}, \cite{li-local}, \cite{ll-book})
that $L(k,0)$ has a natural vertex operator algebra structure
and that all the irreducible $L(k,0)$-modules canonically correspond
to the irreducible highest weight integrable $\hat{\g}$-modules of level $k$.
Let $L(k,\lambda)$ be a highest weight integrable $\hat{\g}$-module 
of level $k$, which is also an $L(k,0)$-module, 
where $\lambda \in \h^*$ is a suitable dominant weight.
We have (cf. \cite{dl})
\begin{equation}
[L(0),a(n)]=-na(n)\;\;\;\mbox{ for }a\in \g,\; n\in \Z,
\end{equation}

\noindent
which implies that $[L(0)-d_{H}, a(n)]=0$ for $a\in \g,\; n\in \Z$.
Consequently, $L(0) = d_H+\mu$ on $L(k,\lambda)$, where 
$\mu$ is the lowest $L(0)$-weight of $L(k,\lambda)$.

Recall the \emph{homogeneous character} of $L(k,\lambda)$ 
\begin{equation}
\chi_{L(k,\lambda)}^H(q) = \tr_{L(k,\lambda)} \, q^{d_H}
\end{equation}

\noindent
and the \emph{principal character} of $L(k,\lambda)$
\begin{eqnarray}
\chi_{L(k,\lambda)}^P(q) =\tr_{L(k,\lambda)} \, q^{-d_P}.
\end{eqnarray}

\bd{character-defs} 
{\em Define a multi-parameter character of $L(k,\lambda)$ 
as follows (cf. \cite{kacbook1}):
\begin{eqnarray}
\chi_{L(k,\lambda)}(x_1,x_2,...,x_{\ell};q) 
= \tr_{L(k,\lambda)} \, x_1^{H^{(1)}(0)}\cdots x_{\ell}^{H^{(\ell)}(0)} q^{L(0)}.
\end{eqnarray}}

\ed

Now, with $L(0) = d_H + \mu$, using Lemma \ref{principal-grading} we have 
\begin{eqnarray}
& &\chi_{L(k,\lambda)}^H(q) = \tr_{L(k,\lambda)} \, q^{d_H} 
= q^{-\mu} \chi_{L(k,\lambda)}(1,\dots,1;q),\\
& &\chi_{L(k,\lambda)}^P(q) = \tr_{L(k,\lambda)} \, q^{-d_P}
= q^{-({\rm ht}(\theta)+1)\mu}
\chi_{L(k,\lambda)}(q^{-1},...,q^{-1};q^{{\rm ht}(\theta)+1}). \label{princhar}
\end{eqnarray}

The following is our main theorem:

\bt{char-recurrence} 
Let $L(k,\lambda)$ be a highest weight integrable
module of level $k$ for an untwisted affine Lie algebra $\hat\g$.
The following relations hold:
\begin{eqnarray} \label{recu-relation}
\chi_{L(k,\lambda)}(x_1,\dots,x_{\ell};q) 
= (x_iq)^{\frac{2}{\<\alpha_i, \alpha_i\>}k}
\chi_{L(k,\lambda)}(x_1q^{a_{1i}},\dots,x_{\ell}q^{a_{\ell i}};q)
\end{eqnarray}

\noindent
for $1 \leq i \leq \ell$, where $C=(a_{ij})$ is the Cartan matrix of 
$\g$. In particular, if $\g$ is simply-laced, we have 
\begin{eqnarray}\label{erecu-relation}
\chi_{L(k,\lambda)}(x_1,\dots,x_{\ell};q) 
= (x_iq)^k \chi_{L(k,\lambda)}(x_1q^{a_{1i}},\dots,x_{\ell}q^{a_{\ell i}};q)
\end{eqnarray}

\noindent
for $1 \leq i \leq \ell$.
\et

\begin{proof} Let $H \in Q^{\vee} \subset \g \subset L(k, 0)$. Set 
\begin{eqnarray}
\Delta(H,x)
=x^{H(0)}\exp \left(\sum_{n\ge 1}\frac{H(n)}{-n}(-x)^{-n}\right).
\end{eqnarray}

\noindent
For $v \in L(k,0)$, set
\begin{equation}
Y_{L(k,\lambda)^{(H)}}(v,x) = Y_{L(k,\lambda)}(\Delta(H,x)v,x) 
= \sum_{n \in \Z} v_{(H)}(n) x^{-n-1}. 
\end{equation}

\noindent
We also set
\begin{equation}
Y_{L(k,\lambda)^{(H)}}(\omega,x) 
= \sum_{n \in \Z} \omega_{(H)}(n) x^{-n-1} 
= \sum_{n \in \Z} L_{(H)}(n) x^{-n-2}.
\end{equation}

It was proved in (\cite{li-affine}, Proposition 2.25) that $(L(k,\lambda), Y_{L(k,\lambda)^{(H)}})$ 
carries the structure of an $L(k,0)$-module and
it (briefly denoted by $L(k,\lambda)^{(H)}$) is isomorphic to $L(k,\lambda)$ as an $L(k,0)$-module.
In particular,
$L(k,\lambda)^{(H_i)}$ is isomorphic to $L(k,\lambda)$ 
for $1 \leq i \leq \ell$. Consequently, we have
\begin{eqnarray}
\chi_{L(k,\lambda)}(x_1,\dots,x_{\ell};q) & = &
\chi_{L(k,\lambda)^{(H_i)}}(x_1,\dots,x_{\ell};q) \notag \\
& = & \tr_{L(k,\lambda)} x_1^{(H^{(1)})_{(H_i)}(0)} 
\dots x_{\ell}^{(H^{(\ell)})_{(H_i)}(0)} q^{L_{(H_i)}(0)}. 
\end{eqnarray}

Next, we find an explicit expression for each of
$(H^{(j)})_{(H_i)}(0)$.
For $G \in \h \subset L(k, 0)$ and $m \in \N$, we have 
\begin{equation}
H(m) G = H(m) G(-1) \1 = G(-1) H(m) \1 + m\<H,G\> \delta_{m,1}k \1.
\end{equation}  
Since $H(m) \1 = 0$ for $m \geq 0$ (from the creation property (\ref{creationprop})), we have
\begin{eqnarray}
H(m) G = \<H,G\> \delta_{m,1} k \1 \;\;\;\mbox{ for } \; m \ge 0,
\end{eqnarray}

\noindent
hence we have $x^{H(0)} G = x^{0} G = G$ and
\begin{eqnarray*}
  \left( \sum_{m \geq 1} \frac{H(m)}{-m} (-x)^{-m} \right) G 
  =\<H,G\> k \1 x^{-1}. 
\end{eqnarray*}

\noindent
Then for $n\ge 2$ we have
\begin{equation}
\left(\sum_{m \geq 1}\frac{H(m)}{-m} (-x)^{-m} \right)^n G = 
k\<H,G\> x^{-1} \left( \sum_{m \geq 1}\frac{H(m)}{-m}(-x)^{-m} \right)^{n-1} {\bf 1} = 0.
\end{equation}

\noindent
Therefore,
\begin{eqnarray*}
\Delta(H,x)G = G + k\<H,G\>x^{-1}{\bf 1},
\end{eqnarray*}

\noindent
{}from which we get
\begin{eqnarray}
 Y_{L(k,\lambda)^{(H)}}(G,x) = Y_{L(k,\lambda)}(\Delta(H,x)G,x) 
= Y_{L(k,\lambda)}(G,x) + \<H,G\> k x^{-1}. 
\end{eqnarray}

\noindent
Taking the residue with respect to $x$ we obtain
\begin{eqnarray}
 G_{(H)}(0) = G(0) + \<H,G\> k. 
\end{eqnarray}

\noindent
In particular, we have
\begin{eqnarray}
(H^{(j)})_{(H_i)}(0) & = & H^{(j)}(0) + \<H_i, H^{(j)}\>k 
 = H^{(j)}(0) + \frac{2}{\<\alpha_i, \alpha_i\>} \delta_{i,j}k.
\end{eqnarray}

\noindent
On the other hand, a similar computation (see \cite{li-affine}, Remark 2.23) gives
\begin{eqnarray}
L_{(H)}(0) = L(0) + H(0) + \frac{\<H,H\>}{2} k.
\end{eqnarray}

\noindent
Thus
\begin{eqnarray}
L_{(H_i)}(0) & = & L(0) + H_i(0) + \frac{\<H_i,H_i\>}{2}k \notag \\
 & = & L(0) + \sum_{j=1}^{\ell} a_{ji} H^{(j)}(0) 
+ \frac{2}{\<\alpha_i, \alpha_i\>}k.
\end{eqnarray}

\noindent
Therefore, we have
\begin{eqnarray*}
\lefteqn{\chi_{L(k,\lambda)}(x_1,\dots,x_{\ell};q)} \\ 
& = & \tr_{L(k,\lambda)} \, x_1^{(H^{(1)})_{(H_i)}(0)} \cdots x_{\ell}^{(H^{(\ell)})_{(H_i)}(0)} q^{L_{(H_i)}(0)} \\ 
& = &\tr_{L(k,\lambda)} \, x_1^{H^{(1)}(0) 
+ \frac{2}{\<\alpha_i, \alpha_i\>} \delta_{i,1}k} 
\cdots x_{\ell}^{H^{(\ell)}(0) 
+ \frac{2}{( \alpha_i, \alpha_i )} \delta_{i,l}k}
q^{L(0)+\sum_{j=1}^{\ell} a_{ji} H^{(j)}(0) 
+ \frac{2}{\<\alpha_i, \alpha_i\>}k} \\
&=& (x_iq)^{\frac{2}{\<\alpha_i, \alpha_i\>}k} 
\tr_{L(k,\lambda)} \, 
(x_1q^{a_{1i}})^{H^{(1)}(0)} \cdots (x_{\ell}q^{a_{\ell i}})^{H^{(\ell)}(0)} q^{L(0)} \\
&=& (x_iq)^{\frac{2}{\<\alpha_i, \alpha_i\>}k} 
\chi_{L(k,\lambda)}(x_1q^{a_{1i}},\dots,x_{\ell}q^{a_{\ell i}};q),
\end{eqnarray*}

\noindent
proving the first assertion. The second assertion follows 
as $\< \alpha_i, \alpha_i \> = 2$ for all 
$1 \leq i \leq \ell$ when  $\g$ is simply-laced.
\end{proof}

Notice that operators $L(0)$ and $H^{1}(0),\dots, H^{\ell}(0)$
on the vertex operator algebra $L(k,0)$ have integral eigenvalues, 
so that
$$\chi_{L(k,0)}(x_1,\dots,x_{\ell};q) \in \C[[x_{1}^{\pm
      1},\dots,x_{\ell}^{\pm 1},q]].$$
For ${\bf n}=(n_{1},\dots,n_{\ell}) \in \Z^{\ell}$,
we set ${\bf x}^{\bf n} = x_1^{n_1} \dots x_{\ell}^{n_{\ell}}$. 
Define $A({\bf n};q) \in \C[[q]]$ by
\begin{eqnarray}
\chi_{L(k,0)}(x_1,\dots,x_{\ell};q) 
&=& \sum_{n_1,\dots,n_{\ell} \in \Z} 
A(n_1,\dots,n_{\ell};q) \, x_1^{n_1}\dots x_{\ell}^{n_{\ell}} \notag \\
& = & \;\;\; \sum_{{\bf n} \in \Z^{\ell}}\;\;\; 
A({\bf n};q) \, {\bf x}^{\bf n}. 
\end{eqnarray}

Then we have:

\bp{pgeneral-recurrsion}

The following recurrence relations hold:
\begin{equation}\label{egeneral-recurrsion}
A(\mathbf{n};q) 
= A\left(n_1,\dots,n_{i-1},n_i-\frac{2}{\<\alpha_i, \alpha_i\>} k,n_{i+1},\dots,n_{\ell};q\right)
q^{- \frac{2}{\<\alpha_i, \alpha_i\>}k + \sum_{m=1}^{\ell} a_{mi} n_m}
\end{equation}
for $\mathbf{n} = (n_1,\dots,n_{\ell}) \in \Z^{\ell}$, $i=1,\dots,\ell$.
In particular, when $\g$ is simply-laced, we have
\begin{eqnarray}\label{e-ADE-recu}
A(\mathbf{n};q) 
= A(n_1,\dots,n_{i-1},n_i-k,n_{i+1},\dots,n_{\ell};q)
q^{-k +\sum_{j=1}^{\ell} a_{ji} n_j}. 
\end{eqnarray}

\ep

\begin{proof}

The recurrence relations (\ref{erecu-relation}) give us:
\begin{eqnarray*}
\lefteqn{ \chi_{L(k,0)}(x_1,\dots,x_{\ell};q) }\\
& = & (x_iq)^{\frac{2}{\<\alpha_i, \alpha_i\>}k} \chi(x_1q^{a_{1i}},\dots,x_{\ell}q^{a_{\ell i}};q) \\
& = & (x_iq)^{\frac{2}{\<\alpha_i, \alpha_i\>}k} \sum_{\mathbf{n} \in \Z^{\ell}} A(\mathbf{n};q)(x_1q^{a_{1i}})^{n_1}\dots (x_{\ell}q^{a_{\ell i}})^{n_{\ell}} \\
& = & \sum_{\mathbf{n} \in \Z^{\ell}} A(\mathbf{n};q)q^{\frac{2}{\<\alpha_i, \alpha_i\>}k + \sum_{m=1}^{\ell} a_{mi} n_m} x_1^{n_1}\dots x_{i-1}^{n_{i-1}} x_i^{n_i + \frac{2}{\<\alpha_i, \alpha_i\>}k} x_{i+1}^{n_{i+1}} \dots x_{\ell}^{n_{\ell}}.
\end{eqnarray*}

\noindent
Since the sum is over all integers and for $1\le i\le \ell$, 
$\frac{2}{\<\alpha_i, \alpha_i\>}$ is a (positive) integer,
replacing $n_i$ with 
$n_i-\frac{2}{\<\alpha_i, \alpha_i\>}k$ we get
\begin{eqnarray*}
\lefteqn{ \chi_{L(k,0)}(x_1,\dots,x_{\ell};q) }\\
& = & \sum_{n_1,\dots,n_{\ell} \in \Z} A\left(n_1,\dots,n_{i-1},n_i
-\frac{2}{\<\alpha_i, \alpha_i\>} k,n_{i+1},\dots,n_{\ell};q\right) \\ 
& & \qquad \qquad \qquad \qquad \qquad \qquad \qquad 
q^{\frac{2}{\<\alpha_i, \alpha_i\>}k 
- a_{ii}\frac{2}{\<\alpha_i, \alpha_i\>}k 
+ \sum_{j=1}^{\ell} a_{ji} n_j} x_1^{n_1} \dots x_{\ell}^{n_{\ell}} \\
& = & \sum_{n_1,\dots,n_{\ell} \in \Z} A\left(n_1,\dots,n_{i-1},n_i
-\frac{2}{\<\alpha_i, \alpha_i\>} k,n_{i+1},\dots,n_{\ell};q\right) \\
& & \qquad \qquad \qquad \qquad \qquad \qquad \qquad 
q^{- \frac{2}{\<\alpha_i, \alpha_i\>}k 
+ \sum_{j=1}^{\ell} a_{ji} n_j} x_1^{n_1} \dots x_{\ell}^{n_{\ell}}.
\end{eqnarray*}

\noindent
Equating the coefficients of ${\bf x}^{\bf n}$ 
we get (\ref{egeneral-recurrsion}).

If $\g$ is simply-laced, we have
$\frac{2}{\<\alpha_i, \alpha_i\>}=1$ for $i=1,\dots,\ell$, 
so that the second assertion follows.
\end{proof}

\br{initial-cond}

{\em In the case that $\g$ is simply-laced, it is clear from 
(\ref{e-ADE-recu}) that 
all the coefficients $A(n_1,\dots,n_{\ell};q)$ can be uniquely 
determined by using $k^{\ell}$ ($q$-series) initial conditions.}

\er

\br{rspecial-case}

{\em Let us consider the special case when $\g$ is simply-laced with
$k=1$. By (\ref{e-ADE-recu}) we have
\begin{eqnarray}
A(n_1,\dots,n_{\ell};q) 
= A(n_1,\dots,n_{i-1},n_i-1,n_{i+1},\dots,n_{\ell};q)
 q^{-1 + \sum_{m=1}^{\ell} a_{mi} n_m}
\end{eqnarray}
for $n_1,\dots,n_{\ell} \in \Z$, $i=1,\dots,\ell$. 
Now we show
\begin{equation}\label{enew}
A(n_1,\dots,n_{\ell};q) = A(0,\dots,0;q) \sum_{{\bf n} \in \Z^{\ell}} 
                          q^{\frac{1}{2} {\bf n} C {\bf n}^{t}}
\end{equation}
for $n_1,\dots,n_{\ell} \in \Z$,
where $C$ is the Cartan matrix of $\g$ as before.
First, we show that for $i=1,\dots,\ell$
\begin{eqnarray} \label{coef-form}
A(n_1,\dots,n_{\ell};q) = A(n_1,\dots,n_{i-1},0,n_{i+1},\dots,n_{\ell};q)
 q^{-n_i^2 + \sum_{j=1}^{\ell} n_j a_{ji} n_i}.
\end{eqnarray}

If $n_i=0$, the formula holds trivially. 
Assume that (\ref{coef-form}) holds for $n_i-1$. Then
we have
\begin{eqnarray*}
\lefteqn{A(n_1,\dots,n_{\ell};q)} \\
& = & A(n_1,\dots,n_{i-1},n_i-1,n_{i+1},\dots,n_{\ell};q)q^{-1 
+ \sum_{j=1}^{\ell} a_{ji} n_j} \\
& = & A(n_1,\dots,n_{i-1},0,n_{i+1},\dots,n_{\ell};q) \\
& &  \qquad \qquad \qquad \qquad \qquad \cdot q^{-(n_i-1)^2 
+ \sum_{j=1}^{\ell} (n_j-\delta_{j,i})a_{ji} (n_i-1)} q^{-1 
+ \sum_{j=1}^{\ell} a_{ji} n_j} \\
& = & A(n_1,\dots,n_{i-1},0,n_{i+1},\dots,n_{\ell};q)
q^{-n_i^2 + \sum_{j=1}^{\ell} n_j a_{ji} n_i}. 
\end{eqnarray*}

In fact, from this it is easy to see that (\ref{coef-form}) holds for 
$n_i-1$ if and only if it holds for $n_i$. 
Hence by induction, (\ref{coef-form}) holds for all $n_i \in \Z$.

Now applying (\ref{coef-form}) repeatedly, we get (\ref{enew}) as
\begin{eqnarray*}
\lefteqn{A(n_1,\dots,n_{\ell};q)} \nonumber\\
& = & A(0,n_2,\dots,n_{\ell};q) q^{-n_1^2 
+ \sum_{j=1}^{\ell} n_j a_{j1} n_1} \nonumber\\
& = & A(0,0,n_3,\dots,n_{\ell};q) q^{-n_1^2-n_2^2 
+ \sum_{j=1}^{\ell} n_j a_{j1} n_1 + \sum_{j=2}^{\ell} n_j a_{j2} n_2} \nonumber \\
& = & \dots \nonumber\\
& = & A(0,\dots,0;q) q^{-\sum_{j=1}^{\ell} n_j^2 
+ \sum_{j=1}^{\ell} n_j a_{j1} n_1 
+ \sum_{j=2}^{\ell} n_j a_{j2} n_2 + \dots + n_1 a_{\ell \ell} n_{\ell}} \nonumber\\
& = & A(0,\dots,0;q) q^{\frac{1}{2} \sum_{i=1}^{\ell} 
\sum_{m_i=1}^{\ell} n_{m_i} a_{m_ii} n_i} \nonumber\\
& = & A(0,\dots,0;q) q^{\frac{1}{2} {\bf n} C {\bf n}^{t}},
\end{eqnarray*}

\noindent
which proves equation (\ref{enew}).}
\er

Recall the explicit construction (see \cite{fk}, \cite{s}, \cite{flm})  
\begin{eqnarray*}
L(1,0)= V_{Q} = S(\hat{\h}^{-}) \otimes \C[Q],
\end{eqnarray*} 
where $\hat{\h}^{-} = \h \otimes t^{-1} \C [t^{-1}]$ 
is an abelian subalgebra of
$\hat{\g}$ and $\C[Q]$ is the the group algebra of the root lattice $Q$.

As $A(0,\dots,0;q)$ is the coefficient of $x_1^0 \dots x_{\ell}^0$
in $\chi_{L(1,0)}(x_1,\dots,x_{\ell};q)$, we have 
$$A(0,\dots,0;q) = \tr_{S(\hat{\h}^{-})} \, 
q^{L(0)}= \prod_{j \geq 1} (1-q^j)^{-\ell}.$$
 
\noindent
Hence by Remark \ref{rspecial-case} we immediately have
\begin{eqnarray}\label{echar-exp}
\chi_{L(1,0)}(x_1,\dots,x_{\ell};q) = 
\prod_{j \geq 1} (1-q^j)^{-\ell} \sum_{{\bf n} \in \Z^{\ell}} 
q^{\frac{1}{2} {\bf n} C {\bf n}^{t}} {\bf x}^{\bf n},
\end{eqnarray}

\noindent
which is well known.

In particular, by (\ref{princhar}) we have
\begin{eqnarray}
\chi_{L(1,0)}^P(q) \label{char-cor}
= \prod_{j \geq 1} (1-q^{({\rm ht}(\theta)+1)j})^{-\ell}
\sum_{{\bf n}\in \Z^{\ell}} q^{\frac{{\rm ht}(\theta)+1}{2} {\bf n} C
{\bf n}^{t} - \sum_{i=1}^{\ell} n_i}.
\end{eqnarray}

\br{multisum-identity}
{\em There is a known product expression of the
principal character of $L(1,0)$ (\cite{L3}; cf. \cite{LM}, \cite{kvstrange}). 
Equating the product form of the principally
specialized character with the multisum expression in (\ref{char-cor})
gives combinatorial identities which are essentially those of Macdonald (\cite{MacD}; cf. \cite{kvstrange}) 
For example, when $\g= A_{\ell}$ we
have the multisum identity:
\begin{equation}\label{identity-a}
 \prod_{j \geq 1} \frac{(1-q^{(\ell+1)j})^{\ell+1}}{(1-q^j)} = 
   \sum_{{\bf n} \in \Z^{\ell}} q^{\frac{\ell+1}{2} {\bf n}C {\bf n}^{t} -
   \sum_{i=1}^{\ell} n_i} 
\end{equation}

\noindent
and when $\g= D_{\ell}$, we have the identity:
\begin{equation}\label{identity-d}
\prod_{j \geq 1} \frac{(1-q^{2(\ell-1)j})^{\ell}}{(1-q^{2j-1})(1-q^{(\ell-1)(2j-1)})}
= \sum_{{\bf n} \in \Z^{\ell}} q^{(\ell-1) {\bf n}C {\bf n}^{t} - \sum_{i=1}^{\ell} n_i}.
\end{equation}

For $\g= E_{6}, E_{7}, E_{8}$ we also have corresponding multisum
identities (see \cite{c-thesis}). 
For $\ell= 1$, identity (\ref{identity-a}) reduces to a well known identity due to Gauss (cf. \cite{A}, page 23,
(2.2.13)), which was interpreted in this way in \cite{fl}, Corollary 4.7, by principal specialization of the 
explicit character of the level $1$ basic representation of the affine Lie algebra $A_1^{(1)}$. In \cite{kvstrange}, 
Equation (3.21), Kac obtained an identity corresponding to each order $m$ ``rational'' automorphism of a finite 
dimensional simple Lie algebra. The above mentioned identity due to Gauss corresponds to an order $4$ rational 
automorphism of the simple Lie algebra $A_1$. It is also interesting to note that
the isomorphism of the two vertex-operator constructions of $\widehat{sl_{2}}$ reflected by this Gauss identity
is found to be a source of ``triality'' in the construction of the moonshine module vertex operator algebra
(see section 4.5 in \cite{flm}).
} \er


\begin{thebibliography}{CLM-bib}

\bibitem[A]{A} 
G. E. Andrews, {\it The theory of partitions}, In:
G.-C. Rota (ed.) Encyclopedia of Mathematics and its Applications,
{\bf Vol. 2}, Reading, MA, Addison-Wesley 1976.

\bibitem[B]{B}
R. E. Borcherds, Vertex algebras, Kac-Moody algebras, and the Monster,
{\em Proc. Natl. Acad. Sci. USA} {\bf 83} (1986), 3068-3071.

\bibitem[Ca]{c}
S. Capparelli, On some representations of twisted affine Lie algebras and 
combinatorial identities, {\em J. Algebra} {\bf 154} (1993), 335-355. 

\bibitem[CLM1]{clm1}
S. Capparelli, J. Lepowsky and A. Milas,
The Rogers-Ramanujan recursion and intertwining operators, 
{\em Commun. Contemp. Math.} {\bf 5} (2003), 947-966.

\bibitem[CLM2]{clm2}
S. Capparelli, J. Lepowsky and A. Milas,
The Rogers-Selberg recursions, the Gordon-Andrews identities and intertwining operators, 
{\em The Ramanujan J.} {\bf 12} (2006), 377-395.
    
\bibitem[Co]{c-thesis}
W. Cook, Affine Lie algebras, vertex operator algebras, and combinatorial identities,
Ph.D. thesis, North Carolina State University, 2005.

\bibitem[DL]{dl}
C. Dong and J. Lepowsky, {\it Generalized Vertex Algebras and Relative Vertex Operators}, 
Progress in Math., {\bf Vol. 112}, Birkh\"auser, Boston, 1993.

\bibitem[DLM]{dlm}
C. Dong, H.-S. Li and G. Mason, 
Regularity of rational vertex operator algebras,
{\em Adv. Math.} {\bf 132} (1997), 148-166.

\bibitem[D]{dyson}
F. J. Dyson, Missed opportunities,
{\em Bull. Amer. Math. Soc.} {\bf 78} (1972), 635-652. 


\bibitem[FL]{fl}
A. J. Feingold and J. Lepowsky, 
The Weyl-Kac character formula and power series identities, 
{\em Adv. Math.} {\bf 29} (1978), 271-309.

\bibitem[FHL]{fhl}
I. Frenkel, Y.-Z. Huang and J. Lepowsky, On axiomatic approaches to
vertex operator algebras and modules, Memoirs Amer. Math. Soc. {\bf 104}, 1993.

\bibitem[FK]{fk}
I. Frenkel and V. G. Kac, Basic representations of affine Lie algebras and dual resonance
models, {\em Invent. Math.} {\bf 62} (1980), 23-66.

\bibitem[FLM]{flm}
I. Frenkel, J. Lepowsky and A. Meurman, {\it Vertex Operator Algebras and the Monster}, 
Pure and Appl. Math., {\bf Vol. 134}, Academic Press, Boston, 1988.

\bibitem[FZ]{fz}
I. Frenkel and Y.-C. Zhu, Vertex operator algebras associated to representations of 
affine and Virasoro algebras, {\em Duke Math. J.} {\bf 66} (1992), 123-168.

\bibitem[GKS]{GKS} 
F. Garvan, D. Kim, and D. Stanton, Cranks and $t$-cores, 
{\em Invent. Math.} {\bf 101} (1990), 1-17.


\bibitem[H]{hum}
J. Humphreys, {\em Introduction to Lie Algebras and Their Representations,}
Springer-Verlag, New York-Heildelberg-Berlin, 1972.

\bibitem[K1]{k1974}
V. G. Kac, Infinite-dimensional Lie algebras and Dedekind's $\eta$-function,
{\em Funk. Anal. i Prilozhen.} {\bf 8} (1974), 77-78; English transl. {\em Funct. Anal. Appl.}
{\bf 8} (1974), 68-70.

\bibitem[K2]{kvstrange}
V. G. Kac, Infinite-dimensional algebras, Dedekind's $\eta$-function, classical
Mobius function and the very strange formula, {\em Adv. Math.} {\bf 30} (1978), 85-136.

\bibitem[K3]{kacbook1}
V. G. Kac, {\it Infinite Dimensional Lie Algebras},
Cambridge University Press, 3rd edition, 1990.

\bibitem[L1]{L1} 
J. Lepowsky, Macdonald-type identities, 
{\em Adv. Math.} {\bf 27} (1978), no. 3, 230-234.

\bibitem[L2]{L2} 
J. Lepowsky, Generalized Verma modules, loop space cohomology and 
Macdonald-type identities, {\em Ann. Sci. Ecole Norm. Sup. (4)} {\bf 12} (1979), no. 2, 169-234.

\bibitem[L3]{L3} 
J. Lepowsky, Application of the numerator formula to $k$-rowed plane partitions,
{\em Adv. Math.} {\bf 35} (1980), no. 2, 179-194.

\bibitem[LL]{ll-book}
J.  Lepowsky and H.-S. Li, {\it Introduction to Vertex Operator
Algebras and Their Representation Theory},  Progress in Math., {\bf Vol. 227}, 
Birkh\"{a}user, Boston, 2004.

\bibitem[LM]{LM}
J. Lepowsky and S. Milne, Lie algebraic approaches to classical partition identities, 
{\em Adv. Math.} {\bf 29} (1978), 15-59.

\bibitem[LP]{LP}
J. Lepowsky and M. Primc, Structure of the standard modules for the affine Lie algebra $A_1^{(1)}$,
{\em Contemporary Math.} {\bf 46}, 1985.

\bibitem[LW1]{LW1}
J. Lepowsky and R. L. Wilson, Construction of the affine Lie algebra $A_1^{(1)}$, 
{\em Commun. Math. Phys.} {\bf 62} (1978), 43-53.

\bibitem[LW2]{LW2}
J. Lepowsky and R. L. Wilson, A new family of algebras underlying the Rogers-Ramanujan 
identities and generalizations, {\em Proc. Natl. Acad. Sci. USA}, {\bf 78} (1981), 7254-7258. 

\bibitem[LW3]{LW3}
J. Lepowsky and R. L. Wilson, A Lie theoretic interpretation and proof of the 
Rogers-Ramanujan identities, {\em Adv. Math.} {\bf 45} (1982), 21-72.

\bibitem[LW4]{LW4}
J. Lepowsky and R. L. Wilson,The structure of standard modules, I, Universal algebras and 
the Rogers-Ramanujan identities, {\em Invent. Math.} {\bf 77} (1984), 199-290.

\bibitem[LW5]{LW5}
J. Lepowsky and R. L. Wilson,The structure of standard modules, II, The case $A_1^{(1)}$, 
principal gradation, {\em Invent. Math.} {\bf 79} (1985), 417-442.


\bibitem[Li1]{li-local}
H.-S. Li, Local systems of vertex operators, vertex superalgebras and
modules, {\em J. Pure Appl. Alg.} {\bf 109} (1996), 143-195.


\bibitem[Li2]{li-affine}
H.-S. Li, Certain extensions of vertex operator algebras of affine type, 
{\em Commun. Math. Phys.} {\bf 217} (2001), 653-696.

\bibitem[Ma]{MacD}
I. G. Macdonald, Affine root systems and Dedekind's $\eta$-function, 
{\em Invent. Math.} {\bf 15} (1972), 91-143.

\bibitem[Mi1]{Mi1}
K. C. Misra, Level one standard modules for affine symplectic Lie algebras, {\em Math. Ann.} {\bf 287} (1990), 287-302.

\bibitem[Mi2]{Mi2}
K. C. Misra, Level two standard $\tilde{A}_n$-modules, {\em J. Algebra} {\bf 137} (1991), 56-76.

\bibitem[MP1]{mpaim}
A. Meurman and M. Primc, Annihilating ideals of standard modules of 
$\widetilde{sl(2,\C)}$ and combinatorial identities,  
{\em Adv. Math.} {\bf 64} (1987), 177-240.

\bibitem[MP2]{mp96}
A. Meurman and M. Primc,  Vertex operator algebras and 
representations of affine Lie algebras, {\em Acta Applicandae Math.}
{\bf 44} (1996), 207-215. 

\bibitem[MP3]{mp98}
A. Meurman and M. Primc
A basis of the basic $\tilde{sl}(3,C)$-module,
{\em Commun. Contemp. Math.} {\bf 3} (2001), 593-614.

\bibitem[MP4]{mpbook}
A. Meurman and M. Primc, Annihilating Fields of Standard Modules of 
$\widetilde{sl(2,\C)}$ and Combinatorial Identities, 
Memoirs Amer. Math. Soc. {\bf 652}, 1999.

\bibitem[Mo]{M} 
R. Moody, Macdonald identities and Euclidean Lie algebras, 
{\em Proc. Amer. Math. Soc.} {\bf 48} (1975), 43-52.

\bibitem[P1]{primc9805}
M. Primc, Some crystal Rogers-Ramanujan type identities,
{\em Glas. Mat. Ser. III} {\bf 34} (1999), 73-86.

\bibitem[P2]{primccombiden}
M. Primc, Vertex algebras and combinatorial identities.  The 2000
Twente Conference on Lie Groups (Enschede), {\em Acta Appl. Math.}
{\bf 73} (2002), 221-238.

\bibitem[S]{s}
G. Segal, Unitary representations of some infinite-dimensional groups, 
{\em Commun. Math. Phys.} {\bf 80} (1981), 301-342.

\end{thebibliography}
\end{document}